\overfullrule=0pt
\centerline {\bf Multiplicity theorems involving functions with non-convex range}\par
\bigskip
\bigskip
\centerline {BIAGIO RICCERI}\par
\bigskip
\bigskip
\centerline {\it Dedicated to the memory of Professor Csaba Varga, with nostalgia}\par
\bigskip
\bigskip
{\bf Abstract.} Here is a sample of the results proved in this paper: Let $f:{\bf R}\to {\bf R}$ be a continuous function, let $\rho>0$ and let $\omega:[0,\rho[\to [0,+\infty[$ be a continuous increasing function such that 
$\lim_{\xi\to \rho^-}\int_0^{\xi}\omega(x)dx=+\infty$. Consider $C^0([0,1])\times C^0([0,1])$ endowed with the
norm
$$\|(\alpha,\beta)\|=\int_0^1|\alpha(t)|dt+\int_0^1|\beta(t)|dt\ .$$
Then, the following assertions are equivalent:\par
\noindent
$(a)$\hskip 5pt the restriction of $f$ to $\left [-{{\sqrt{\rho}}\over {2}},{{\sqrt{\rho}}\over {2}}\right ]$ is not constant\ ;\par
\noindent
$(b)$\hskip 5pt for every convex set $S\subseteq C^0([0,1])\times C^0([0,1])$ dense in $C^0([0,1])\times C^0([0,1])$,
 there exists $(\alpha,\beta)\in S$ such that the problem
$$\cases{-\omega\left(\int_0^1|u'(t)|^2dt\right)u''=\beta(t)f(u)+\alpha(t) & in $[0,1]$\cr & \cr u(0)=u(1)=0\cr & \cr
\int_0^1|u'(t)|^2dt<\rho\cr}$$
has at least two classical solutions.\par
\bigskip
\bigskip
\bigskip
\bigskip
{\bf 1. - Introduction}\par
\bigskip
Let $H$ be a real Hilbert space. A very classical result of Efimov and Stechkin ([3]) states that if $X$ is a non-convex sequentially weakly
closed subset of $H$, then there exists $y_0\in H$ such that the restriction to $X$ of the function $x\to \|x-y_0\|$ has at least two global
minima. A more precise version of such a result was obtained by I. G. Tsar'kov in [10]. Actually, he proved that any convex set dense in $H$
contains a point $y_0$ with the above property.\par
\smallskip
In the present paper, as a by product of a more general result, we get the following:\par
\medskip
THEOREM 1.1 - {\it Let $X\subset H$ be a non-convex sequentially weakly closed set and let $u_0\in \hbox {\rm conv}(X)\setminus X$.\par
Then, if we put
$$\delta:=\hbox {\rm dist}(u_0,X)$$
and, for each $r>0$,
$$\rho_r:=\sup_{\|y\|<r}((\hbox {\rm dist}(u_0+y,X))^2-\|y\|^2)\ ,$$
for every convex set $S\subseteq H$ dense in $H$, for every bounded sequentially weakly lower semicontinuous function $\varphi:X\to {\bf R}$ and
for every $r$ satisfying
$$r>{{\rho_r-\delta^2+\sup_X\varphi-\inf_X\varphi}\over {2\delta}}\ ,$$
there exists $y_0\in S$, with $\|y_0-u_0\|<r$, such that the function $x\to \|x-y_0\|^2+\varphi(x)$ has at least two global minima
in $X$.}
\medskip
So, with respect to the Efimov-Stechkin-Tsar'kov result, Theorem 1.1 gives us two remarkable 
additional informations: a precise localization of the point
$y_0$ and the validity of the conclusion not only for the function $x\to \|x-y_0\|^2$, but also for suitable perturbations of it.\par
\smallskip
Let us recall the most famous open problem in this area: if $X$ is a subset of $H$ such that, for
each $y\in H$, the restriction of the function $x\to \|x-y\|$ to $X$ has a unique global minimum, is it true that the set $X$ is convex ? So, Efimov-Stechkin's
result provides an affirmative answer when $X$ is sequentially weakly closed. However, it is a quite common feeling that the answer, in general,
should be negative ([1], [2], [5], [8]). In the light of Theorem 1.1, we posit the following problem:\par
\medskip
PROBLEM 1.1. - Let $X$ be a subset of $H$ for which there exists a bounded sequentially weakly lower semicontinuous function
$\varphi:X\to {\bf R}$ such that, for each $y\in H$, the function $x\to \|x-y\|^2+\varphi(x)$ has a unique global minimum in $X$. Then,
must $X$ be convex ?\par
\medskip
 What allows us to reach the advances presented in Theorem 1.1 is our particular approach which is entirely based on the minimax theorem established in [9]. So, also the present paper can be regarded as a further ring of the chain of applications and consequences of that minimax theorem.\par
\bigskip
{\bf 2. - Results}\par
\bigskip
In the sequel, $X$ is a topological space and $E$ is real normed space, with topological dual $E^*$. 
\smallskip
For each $S\subseteq E^*$, we denote by ${\cal A}(X,S)$ (resp. ${\cal A}_s(X,S))$ 
 the class of all pairs $(I,\psi)$, whith $I:X\to {\bf R}$ and
$\psi:X\to E$, such that, for each $\eta\in S$ and each $s\in {\bf R}$, the set 
$$\{x\in X : I(x)+\eta(\psi(x))\leq s\}$$ is closed and compact (resp. sequentially closed and sequentially compact).\par
\smallskip
Let us start establishing the following useful proposition. $E^{'}$ denotes the algebraic dual of $E$.\par
\medskip
PROPOSITION 2.1. - {\it Let $I:X\to {\bf R}$, let $\psi:X\to E$ and let $x_1,...,x_n\in X$, $\lambda_1,...,\lambda_n\in [0,1]$, with $\sum_{i=1}^n\lambda_i=1$. \par
Then, one has
$$\sup_{\eta\in E^{'}}\inf_{x\in X}\left (I(x)+\eta\left (\psi(x)-\sum_{i=1}^n\lambda_i\psi(x_i)\right )\right )\leq\max_{1\leq i\leq n}I(x_i)\ .$$}\par
\smallskip
PROOF.  Fix $\eta\in E^{'}$. 
Clearly, for some $j'\in \{1,...,n\}$, we have
$$\eta\left (\psi(x_{j'})-\sum_{i=1}^n\lambda_i\psi(x_i)\right )\leq 0\ .\eqno{(1)}$$
Indeed, if not, we would have
$$\eta(\psi(x_j))>\sum_{i=1}^n\lambda_i\eta(\psi(x_i))$$
for each $j\in \{1,...,n\}$. So, multiplying by $\lambda_j$ and summing, we would obtain 
$$\sum_{j=1}^n\lambda_j\eta(\psi(x_j))>\sum_{i=1}^n\lambda_i\eta(\psi(x_i))\ ,$$
a contradiction. 
In view of $(1)$, we have
$$\inf_{x\in X}\left (I(x)+\eta\left (\psi(x)-\sum_{i=1}^n\lambda_i\psi(x_i)\right )\right )
\leq I(x_{j'})+\eta\left (\psi(x_{j'})-\sum_{i=1}^n\lambda_i\psi(x_i)\right )\leq I(x_{j'})\leq \max_{1\leq i\leq n}I(x_i)$$
and so we get the conclusion due to the arbitrariness of $\eta$.\hfill $\bigtriangleup$\par
\medskip
Our main result is as follows:\par
\medskip
THEOREM 2.1. - {\it Let $I:X\to {\bf R}$, let $\psi:X\to E$, let $S\subseteq E^*$ be a convex set
dense in $E^*$ and let $u_0\in E$.\par
Then, for every bounded 
 function $\varphi:X\to {\bf R}$ such that $(I+\varphi,\psi)\in {\cal A}(X,S)$ and for every $r$
satisfying
$$\sup_X\varphi-\inf_X\varphi<\inf_{x\in X}\left( I(x)+\|\psi(x)-u_0\|r\right )-\sup_{\|\eta\|_{E^*}<r}\inf_{x\in X}(I(x)+\eta(\psi(x)-u_0))\ ,
  \eqno{(2)}$$
there exists $\tilde\eta\in S$, with $\|\tilde\eta\|_{E^*}<r$, such that the
function $I+\tilde\eta\circ\psi+\varphi$ has at least two global minima in $X$.}\par
\smallskip
PROOF. Consider the function $g:X\times E^*\to {\bf R}$ defined by
$$g(x,\eta)=I(x)+\eta(\psi(x)-u_0)$$
for all $(x,\eta)\in X\times E^*$. Let $B_r$ denote the open ball in $E^*$, of
radius $r$, centered at $0$. Clearly, for each $x\in X$,
we have
$$\sup_{\eta\in B_r}\eta(\psi(x)-u_0))=\|\psi(x)-u_0\|r\ .\eqno{(3)}$$
Then, from $(2)$ and $(3)$, it follows
$$\sup_X\varphi-\inf_X\varphi<\inf_X\sup_{B_r}g-\sup_{B_r}\inf_Xg\ .\eqno{(4)}$$
Now, consider the function $f:X\times (S\cap B_r)\to {\bf R}$ defined by
$$f(x,\eta)=g(x,\eta)+\varphi(x)$$
for all $(x,\eta)\in X\times (S\cap B_r)$.
Since $S$ is dense in $E^*$, the set $S\cap B_r$ is dense in $B_r$. Hence, 
since $g(x,\cdot)$ is continuous, we obtain
$$\inf_X\sup_{S\cap B_r}g=\inf_X\sup_{B_r}g\ .\eqno{(5)}$$
Then, taking $(4)$ and $(5)$ into account, we have
$$\sup_{S\cap B_r}\inf_Xf\leq \sup_{B_r}\inf_Xf\leq \sup_{B_r}\inf_Xg+\sup_X\varphi<\inf_X\sup_{B_r}g+\inf_X\varphi\leq 
\inf_{x\in X}\left (\sup_{\eta\in S\cap B_r}g(x,\eta)+\varphi(x)\right )=\inf_X\sup_{S\cap B_r}f\ .\eqno{(6)}$$
Now, since $(I+\varphi,\psi)\in {\cal A}(X,S)$ and $f$ is concave in $S\cap B_r$,
we can apply Theorem 1.1 of [9]. Therefore,
since (by $(6)$) $\sup_{S\cap B_r}\inf_Xf<\inf_X\sup_{S\cap B_r}f$, there exists
of $\tilde\eta\in S\cap B_r$ such that the function $f(\cdot,\tilde\eta)$ has at least two global minima in $X$ which, of course, are
global minima of the function $I+\tilde\eta\circ\psi+\varphi$.
\hfill $\bigtriangleup$\par
\medskip
If we renounce to the very detailed informations contained in its conclusion, we can state Theorem 2.1 in an extremely simplified form. 
\medskip
THEOREM 2.2. - {\it Let $I:X\to {\bf R}$, let $\psi:X\to E$ and let $S\subset E^*$ be a convex set 
weakly-star dense in $E^*$. Assume that $\psi(X)$ is not convex and that $(I,\psi)\in {\cal A}(X,S)$.\par
Then, there exists $\tilde \eta\in S$
such that the function $I+\tilde\eta\circ\psi$ has at least two global minima in $X$.}
\smallskip
PROOF. Fix $u_0\in \hbox {\rm conv}(\psi(X))\setminus \psi(X)$ and consider the function $g:X\times E^*\to {\bf R}$ defined by
$$g(x,\eta)=I(x)+\eta(\psi(x)-u_0)$$
for all $(x,\eta)\in X\times E^*$. By Proposition 2.1, we know that
$$\sup_{E^*}\inf_Xg<+\infty\ .$$
On the other hand, for each $x\in X$, since $\psi(x)\neq u_0$, we have
$$\sup_{\eta\in E^*}\eta(\psi(x)-u_0)=+\infty\ .$$
Hence, since $S$ is weakly-star dense in $E^*$ and $g(x,\cdot)$ is weakly-star continuous, we have
$$\sup_{\eta\in S}g(x,\eta)=+\infty\ .$$
Therefore
$$\sup_{S}\inf_Xg<\inf_X\sup_Sg\ .\eqno{(7)}$$
Now, taken into account that $(I,\psi)\in {\cal A}(X,S)$, we can apply Theorem 1.1 of [9] to $g_{|_{X\times S}}$. So,
in view of $(7)$, there exists $\tilde\eta\in S$ such that the function $g(\cdot,\tilde\eta)$ (and so $I+\tilde\eta\circ \psi$) has
at least two global minima in $X$, as claimed.\hfill $\bigtriangleup$\par
\medskip
The next result is a sequential version of Theorem 1.1 of [9].\par
\medskip
THEOREM 2.3. - {\it Let $X$ be a topological space, $E$ a topological vector
space, $Y\subseteq E$ a non-empty separable convex set 
and $f:X\times Y\to
{\bf R}$ a function satisfying the following conditions:\par
\noindent
$(a)$\hskip 5pt for each $y\in Y$, the function $f(\cdot,y)$ is sequentially lower semicontinuous, sequentially
inf-compact and has a unique global minimum in $X$\ ;\par
\noindent
$(b)$\hskip 5pt for each $x\in X$, the function $f(x,\cdot)$ is continuous and quasi-concave.\par
Then, one has
$$\sup_Y\inf_Xf=\inf_X\sup_Yf\ .$$}\par
\smallskip
PROOF. The pattern of the proof is the same as that of Theorem 1.1 of [9]. We limit ourselves to stress the needed changes. First, for every
$n\in {\bf N}$, one
proves the result when $E={\bf R}^n$ and $Y=S_n:=\{(\lambda_1,...,\lambda_n)\in ([0,+\infty[)^n : \lambda_1+...+\lambda_n=1\}$. In this connection,
the proof agrees exactly with that of Lemma 2.1 of [9], with the only difference of using the sequential version of Theorem 1.A of [9] instead
of such a result itself (see Remark 2.1 of [9]). Next, we fix a sequence $\{x_n\}$ dense in $Y$. For each $n\in {\bf N}$, set
$$P_n=\hbox {\rm conv}(\{x_1,...,x_n\})\ .$$
Consider the function $\eta : S_n\to P$ defined by
$$\eta(\lambda_1,...,\lambda_n)=\lambda_1x_1+...+\lambda_nx_n$$
for all $(\lambda_1,...,\lambda_n)\in S_n$. Plainly, the  function
$(x,\lambda_1,...,\lambda_n)\to f(x,\eta(\lambda_1,...,\lambda_n))$ satisfies in $X\times S_n$
the assumptions of Theorem A, and so, by the case previously proved, we have
$$\sup_{(\lambda_1,...,\lambda_n)\in S_n}\inf_{x\in X}f(x,\eta(\lambda_1,...,\lambda_n))=
\inf_{x\in X}\sup_{(\lambda_1,...,\lambda_n)\in S_n}f(x,\eta(\lambda_1,...,\lambda_n))\ .$$
Since $\eta(S_n)=P_n$, we then have
$$\sup_{P_n}\inf_Xf=\inf_X\sup_{P_n}f\ .$$
Now, set
$$D=\bigcup_{n\in {\bf N}}P_n\ .$$
In view of Proposition 2.2 of [9], we have
$$\sup_D\inf_Xf=\inf_X\sup_Df\ .$$
Finally, by continuity and density, we have
$$\sup_{y\in D}f(x,y)=\sup_{y\in Y}f(x,y)$$
for all $x\in X$, and so
$$\inf_X\sup_Yf=\inf_X\sup_Df=\sup_D\inf_Xf\leq\sup_Y\inf_Xf\leq\inf_X\sup_Yf$$
and the proof is complete.\hfill $\bigtriangleup$\par
\medskip
Reasoning as in the proof of Theorem 2.1 and using Theorem 2.3, we get
\medskip
THEOREM 2.4. - {\it Let the assumptions of Theorem 2.1 be satisfied. In addition, assume that
$E^*$ is separable.\par
Then, the conclusion of Theorem 2.1 holds with ${\cal A}_s(X,S)$ instead of ${\cal A}(X,S)$}.\par
\medskip
Analogously, the sequential version of Theorem 2.2 is as follows:\par
\medskip
THEOREM 2.5. - {\it  Let $I:X\to {\bf R}$, let $\psi:X\to E$ and let $S\subset E^*$ be a convex set 
weakly-star separable and weakly-star dense in $E^*$. Assume that $\psi(X)$ is not convex and that $(I,\psi)\in {\cal A}_s(X,S)$.\par
Then, there exists $\tilde \eta\in S$
such that the function $I+\tilde\eta\circ\psi$ has at least two global minima in $X$.}\par
\medskip
Here is a consequence of Theorem 2.1:\par
\medskip
THEOREM 2.6. - {\it Let $E$ be a Hilbert space, let $\psi:X\to E$ be a weakly continuous function and let $S\subseteq E$ be a convex set
dense in $E$. Assume that $\psi(X)$ is not convex and that the function $\|\psi(\cdot)\|$ is inf-compact. Let $u_0\in \hbox {\rm conv}(\psi(X))\setminus \psi(X)$.\par
Then, 
for every bounded 
 function $\varphi:X\to {\bf R}$ such that $\|\psi(\cdot)\|^2+\varphi(\cdot)$ is lower semicontinuous and for every $r$
satisfying
$$r>{{\sup_{\|y\|<r}((\hbox {\rm dist}(u_0+y,\psi(X)))^2-\|y\|^2)-(\hbox {\rm dist}(u_0,\psi(X)))^2+\sup_X\varphi-\inf_X\varphi}\over
{2\hbox {\rm dist}(u_0,\psi(X))}}\ ,\eqno{(8)}$$
there exists $\tilde y\in S$, with $\|\tilde y-u_0\|<r$, such that the function $\|\psi(\cdot)-\tilde y\|^2+\varphi(\cdot)$ has at least two global
minima in $X$.}\par
\smallskip
PROOF. First, we observe that the set $\psi(X)$ is sequentially weakly closed (and so norm closed). Indeed, let $\{x_n\}$ be a sequence
in $X$ such that $\{\psi(x_n)\}$ converges weakly to $y\in E$. So, in particular, $\{\psi(x_n)\}$ is bounded and hence, since $\|\psi(\cdot)\|$
is inf-compact, there exists a compact set $K\subseteq X$ such that $x_n\in K$ for all $n\in {\bf N}$. Since $\psi$ is weakly continuous, the set
$\psi(K)$ is weakly compact and hence weakly closed. Therefore, $y\in \psi(K)$, as claimed. This remark ensures that $\hbox {\rm dist}(u_0,\psi(X))>0$.
Now, we apply Theorem 2.1 identifying $E$ with $E^*$ and taking
$$I(x)={{1}\over {2}}\|\psi(x)-u_0\|^2$$
for all $x\in X$. Of course, we have
$$I(x)+\langle \psi(x)-u_0,y\rangle={{1}\over {2}}(\|\psi(x)-u_0+y\|^2-\|y\|^2) \eqno{(9)}$$
for all $y\in E$.
In view of $(8)$ and $(9)$, we have
$${{1}\over {2}}(\sup_X\varphi-\inf_X\varphi)<{{1}\over {2}}(\hbox {\rm dist}(u_0,\psi(X)))^2+r\hbox {\rm dist}(u_0,\psi(X))-
{{1}\over {2}}\sup_{\|y\|<r}((\hbox {\rm dist}(u_0-y,\psi(X)))^2-\|y\|^2)$$
$$\leq\inf_{x\in X}(I(x)+\|\psi(x)-u_0\|r)-{{1}\over {2}}\sup_{\|y\|<r}((\hbox {\rm dist}(u_0-y,\psi(X)))^2-\|y\|^2)$$
$$=\inf_{x\in X}(I(x)+\|\psi(x)-u_0\|r)-\sup_{\|y\|<r}\inf_{x\in X}(I(x)+\langle \psi(x)-u_0,y\rangle)\ .\eqno{(10)}$$
Let us show that $(I+{{1}\over {2}}\varphi,\psi)\in {\cal A}(X,E)$. So, fix $y\in E$. Since $\psi$ is weakly continuous, $\langle \psi(\cdot),v\rangle$
is continuous in $X$ for all $v\in E$. Observing that
$$I(x)+{{1}\over {2}}\varphi(x)+\langle\psi(x),y\rangle={{1}\over {2}}(\|\psi(x)\|^2+\varphi(x))+\langle\psi(x),y-u_0\rangle+{{1}\over {2}}
\|u_0\|^2\ ,$$
we infer that $I(\cdot)+{{1}\over {2}}\varphi(\cdot)+\langle\psi(\cdot),y\rangle$ is lower semicontinuous since $\|\psi(\cdot)\|^2+\varphi(\cdot)$
is so by assumption. Now, let $s\in {\bf R}$. We readily have
$$\left\{x\in X : I(x)+{{1}\over {2}}\varphi(x)+\langle\psi(x),y\rangle\leq s\right\}\subseteq
\left\{x\in X : \|\psi(x)\|^2-2\|y-u_0\|\|\psi(x)\|\leq 2s-\inf_X\varphi\right\}\ .\eqno{(11)}$$
Since $\|\psi(\cdot)\|$ is inf-compact, the set in the right-hand side of $(11)$ is compact and hence so is the set in left-hand right, as
claimed. Since the set $u_0-S$ is convex and dense in $E$, in view of $(10)$, Theorem 2.1 ensures the existence of $\tilde v\in u_0-S$,
with $\|\tilde v\|<r$, such
that the function $I(\cdot)+\langle\psi(\cdot),\tilde v\rangle+{{1}\over {2}}\varphi(\cdot)$ has at least two global minima in $X$.
Consequently, since
$$I(x)+\langle\psi(x),\tilde v\rangle+{{1}\over {2}}\varphi(x)={{1}\over {2}}(\|\psi(x)+\tilde v-u_0\|^2+\varphi(x))-{{1}\over {2}}
(\|u_0\|^2-\|\tilde v-u_0\|^2)\ ,$$
if we put
$$\tilde y:=u_0-\tilde v\ ,$$
we have $\tilde y\in S$, $\|\tilde y-u_0\|<r$ and the function $\|\psi(\cdot)-\tilde y\|^2+\varphi(\cdot)$ has at least two global minima in $X$.
The proof is complete.\hfill $\bigtriangleup$\par
\medskip
REMARK 2.1. - Of course, Theorem 1.1 is an immediate corollary of Theorem 2.6: take $E=H$, consider $X$ equipped with the relative
weak topology, take $\psi(x)=x$ and observe that if $\varphi:X\to {\bf R}$ is sequentially weakly lower semicontinuous, then $\|\cdot\|^2+
\varphi(\cdot)$ is weakly lower semicontinuous in view of the Eberlein-Smulyan theorem.\par
\medskip
Here is an application of Theorem 2.2. An operator $T$ between two Banach spaces $F_1, F_2$ is said to be sequentially weakly continuous if,
for every sequence $\{x_n\}$ in $F_1$ weakly convergent to $x\in F_1$, the sequence $\{T(x_n)\}$ converges weakly to $T(x)$ in $F_2$.\par
\medskip
THEOREM 2.7. - {\it Let $V$ be a reflexive real Banach space, let $x_0\in V$, let $r>0$, let $X$ be the open ball in $V$, of
radius $r$, centered at $x_0$, let $\gamma:[0,r[\to {\bf R}$, with $\lim_{\xi\to r^-}\gamma(\xi)=+\infty$, let $I:X\to {\bf R}$ and $\psi:X\to E$ be two G\^ateaux differentiable functions. Moreover, assume that $I$ is sequentially weakly lower semicontinous, that $\psi$ is sequentially weakly continuous,
that $\psi(X)$ is bounded and non-convex, 
and that
$$\gamma(\|x-x_0\|)\leq I(x)$$
for all $x\in X$.\par
Then,  for every convex set $S\subseteq E^*$ 
weakly-star dense in $E^*$, there exists $\tilde \eta\in S$ such that
the equation
$$I'(x)+(\tilde\eta\circ\psi)'(x)=0$$
has at least two solutions in $X$.}\par
\smallskip
PROOF. We apply Theorem 2.2 considering $X$ equipped with the relative weak topology. Let $\eta\in E^*$. Since $\psi(X)$ is bounded,
we have $c:=\inf_{x\in X}\eta(\psi(x))>-\infty\ .$ Let $s\in {\bf R}$. We have
$$\{x\in X : I(x)+\eta(\psi(x))\leq s\}\subseteq \{x\in X : I(x)\leq s-c\}\subseteq \{x\in X: \gamma(\|x-x_0\|)\leq s-c\}\ .\eqno{(12)}$$
Since $\lim_{\xi\to r^-}\gamma(t)=+\infty$, there is $\delta\in ]0,r[$, such that $\gamma(\xi)>s-c$ for all $\xi\in ]\delta,r[$. 
Consequently,
from $(12)$, we obtain
$$\{x\in X : I(x)+\eta(\psi(x))\leq s\}\subseteq \{x\in V : \|x-x_0\|\leq \delta\}\ .\eqno{(13)}$$
From the assumptions, it follows that the function $I+\eta\circ\psi$ is sequentially weakly lower semicontinuous in $X$. Hence, from $(13)$, since
$\delta<r$ and $V$ is reflexive, we infer that the set $\{x\in X : I(x)+\eta(\psi(x))\leq s\}$ is sequentially weakly compact and hence weakly compact,
by the Eberlein-Smulyan theorem. In other words, $(I,\psi)\in {\cal A}(X,E^*)$. Therefore, we can apply Theorem 2.2. Accordingly,
there exists $\tilde\eta\in S$ such that the function $I+\tilde\eta\circ\psi$ has at least two global minima in $X$ which are critical points of it since
$X$ is open.\hfill$\bigtriangleup$\par
\medskip
 Here is an application of Theorem 1.1:\par
\medskip
THEOREM 2.8. - {\it Let $H$ be a Hilbert space and let $I, J:H\to {\bf R}$ be two $C^1$ functionals with compact derivative such that $2I-J^2$ is bounded. Moreover, assume that $J(0)\neq 0$ and that there is $\hat x\in H$ such that $J(-\hat x)=-J(\hat x)$.\par
Then, for every convex set $S\subseteq H\times {\bf R}$ dense in $H\times {\bf R}$ and for every $r$ satisfying
$$r>{{\|\hat x\|^2+|J(\hat x)|^2-\inf_{x\in H}(\|x\|^2+|J(x)|^2)+\sup_H(2I-J^2)-\inf_X(2I-J^2)}\over {2\inf_{x\in H}\sqrt {\|x\|^2+|J(x)|^2}}}
\ ,$$
there exists $(y_0,\mu_0)\in S$, with $\|y_0\|^2+|\mu_0|^2<r^2$, such that
the equation
$$x+I'(x)+\mu_0 J'(x) = y_0$$
has at least three solutions.}\par
\smallskip
PROOF. We consider the Hilbert space $E:=H\times {\bf R}$ with the scalar product
$$\langle (x,\lambda),(y,\mu)\rangle_E = \langle x,y\rangle +\lambda\mu$$
for all $(x,\lambda), (y,\mu)\in E$. Take
$$X=\{(x,\lambda)\in E : \lambda=J(x)\}\ .$$
Since $J'$ is compact, the functional $J$ turns out to be sequentially weakly continuous ([11], Corollary 41.9). So, the set $X$ is sequentially weakly closed. Moreover, notice that $(0,0)\not\in X$, while the antipodal points $(\hat x,J(\hat x))$ and $-(\hat x,J(\hat x))$ lie
in $X$. So, $(0,0)\in \hbox {\rm conv}(X)$. Now, with the notations of Theorem 1.1, taking, of course, $u_0=(0,0)$, we have
$$\delta=\inf_{x\in X}\sqrt{\|x\|^2+|J(x)|^2}$$
and
$$\rho_r=\sup_{\|y\|^2+|\mu|^2<r^2}\inf_{x\in X}(\|x\|^2+|J(x)|^2-2\langle (x,J(x)),(y,\mu)\rangle_E)\ .$$
Then, from Proposition 2.1, we infer that
$$\rho_r\leq \|\hat x\|^2+|J(\hat x)|^2\ .$$
Now, consider the function $\varphi:X\to {\bf R}$ defined by
$$\varphi(x,\lambda)=2I(x)-\lambda^2$$
for all $(x,\lambda)\in X$. Notice that $\varphi$ is sequentially weakly continuous and
$r$ satisfies the inequality of Theorem 1.1. Consequently, there exists $(y_0,\mu_0)\in S$ such that the 
functional $$(x,\lambda)\to \|(x,\lambda)\|_E^2-2\langle
(x,\lambda),(y_0,\mu_0)\rangle_E+2I(x)-\lambda^2$$ has at least two global minima in $X$. Of course, if $(x,\lambda)\in X$, we have
$$\|(x,\lambda)\|_E^2-2\langle (x,\lambda),(y_0,\mu_0)\rangle_E+2I(x)-\lambda^2=\|x\|^2+J^2(x)-2\langle x,y_0\rangle-
2\mu_0J(x)+2I(x)-J^2(x)\ .$$
In other words, the functional $x\to \|x\|^2-2\langle x,y_0\rangle-2\mu_0J(x)+2I(x)$ has two global minima in $H$. Since
the functional $x\to -2\langle x,y_0\rangle-2\mu_0J(x)+2I(x)$ has a compact derivative, a well know result ([11], Example 38.25) ensures that
the functional $x\to \|x\|^2-2\langle x,y_0\rangle-2\mu_0J(x)+2I(x)$ has the Palais-Smale property and so, by Corollary 1 of [6], it possesses
at least three critical points. The proof is complete.\hfill $\bigtriangleup$\par
\medskip
REMARK 2.2. - In Theorem 2.8, apart from being $C^1$ with compact derivative, the truly essential assumption on $J$ is, of course, that its graph is not convex. This amounts to say that $J$ is not affine. The current assumptions are made to simplify the constants appearing in the
conclusion. Actually, from the proof of Theorem 2.8, the following can be obtained:\par
\medskip
THEOREM 2.9. - {\it Let $H$ be a Hilbert space and let $I, J:H\to {\bf R}$ be two $C^1$ functionals with compact derivative such that $2I-J^2$ is bounded. Moreover, assume
that $J$ is not affine.\par
Then, for every convex set $S\subseteq H\times {\bf R}$ dense in $H\times {\bf R}$, there exists $(y_0,\lambda_0)\in S$ such that
the equation
$$x+I'(x)+\lambda_0 J'(x) = y_0$$
has at least three solutions.}\par
\medskip
REMARK 2.3. - For $I=0$, the conclusion of Theorem 2.9 can be obtained from Theorem 4 of [7] (see also [4]) provided that, for some $r\in {\bf R}$,
the set $J^{-1}(r)$ is not convex. Therefore, for instance, the fact that, for any non-constant
bounded $C^1$ function $J:{\bf R}\to {\bf R}$, there are $a, b\in {\bf R}$ such that the equation
$$x+aJ'(x)=b$$
has at least three solutions, follows, in any case, from Theorem 2.9, while it follows from Theorem 4 of [7] only if $J$ is not monotone.\par
\medskip
We conclude presenting an application of Theorem 2.7 to a class of Kirchhoff-type problems.\par
\medskip
THEOREM 2.10. - {\it Let $f:{\bf R}\to {\bf R}$ be a continuous function, let $\rho>0$ and let $\omega:[0,\rho[\to [0,+\infty[$ be a continuous increasing function such that $\lim_{\xi\to \rho^-}\int_0^{\xi}\omega(x)dx=+\infty$. Consider $C^0([0,1])\times C^0([0,1])$ endowed with the
norm
$$\|(\alpha,\beta)\|=\int_0^1|\alpha(t)|dt+\int_0^1|\beta(t)|dt\ .$$
Then, the following assertions are equivalent:\par
\noindent
$(a)$\hskip 5pt the restriction of $f$ to $\left [-{{\sqrt{\rho}}\over {2}},{{\sqrt{\rho}}\over {2}}\right ]$ is not constant\ ;\par
\noindent
$(b)$\hskip 5pt for every convex set $S\subseteq C^0([0,1])\times C^0([0,1])$ dense in $C^0([0,1])\times C^0([0,1])$,
 there exists $(\alpha,\beta)\in S$ such that the problem
$$\cases{-\omega\left(\int_0^1|u'(t)|^2dt\right)u''=\beta(t)f(u)+\alpha(t) & in $[0,1]$\cr & \cr u(0)=u(1)=0\cr & \cr
\int_0^1|u'(t)|^2dt<\rho\cr}$$
has at least two classical solutions.}\par
\smallskip
PROOF. Consider the Sobolev space $H^1_0(]0,1[)$ with the usual scalar product
$$\langle u,v\rangle=\int_0^1u'(t)v'(t)dt\ .$$
Let $B_{\sqrt{\rho}}$ be the open ball in $H^1_0(]0,1[$, of radius $\sqrt{\rho}$, centered at $0$ .  Let
$g:[0,1]\times {\bf R}\to {\bf R}$ be a continuous function. Consider the functionals $I, J_g:B_{\sqrt{\rho}}\to {\bf R}$
defined by
$$I(u)={{1}\over {2}}\tilde\omega\left(\int_0^1|u'(t)|^2dt\right)\ ,$$
$$J_g(u)=\int_0^1\tilde g(t,u(t))dt$$
for all $u\in B_{\sqrt{\rho}}$, where $\tilde\omega(\xi)=\int_0^{\xi}\omega(x)dx$, $\tilde g(t,\xi)=\int_0^{\xi}g(t,x)dx$. By classical
results, taking into account that if $\omega(x)=0$ then $x=0$, it follows that the classical solutions of the problem
$$\cases{-\omega\left(\int_0^1|u'(t)|^2dt\right)u''=g(t,u) & in $[0,1]$\cr & \cr u(0)=u(1)=0\cr & \cr
\int_0^1|u'(t)|^2dt<\rho\cr}$$
are exactly the critical points in $B_{\sqrt{\rho}}$ of the functional $I-J_g$. \par
Let us prove that $(a)\to (b)$. 
We are going to apply Theorem 2.7 taking $V=H^1_0(]0,1[)$, $x_0=0$, $r=\sqrt{\rho}$,
$I$ as above, $\gamma(\xi)={{1}\over {2}}\tilde\omega(\xi^2)$,
$E=C^0([0,1])\times C^0([0,1])$ and
$\psi:B_{\sqrt{\rho}}\to E$ defined by 
$$\psi(u)(\cdot)=(u(\cdot),\tilde f(u(\cdot)))$$
for all $u\in B_{\sqrt{\rho}}$, where $\tilde f(\xi)=\int_0^{\xi}f(x)dx$. Clearly, the functional $I$ is continuous and strictly convex (and so weakly lower semicontinuous), while the operator $\psi$ is G\^ateaux differentiable  and sequentially weakly continuous due to the compact embedding of $H^1_0(]0,1[)$
into $C^0([0,1])$. Recall that
$$\max_{[0,1]}|u|\leq {{1}\over {2}}\sqrt{\int_0^1|u'(t)|^2dt}$$
for all $u\in H^1_0(]0,1[)$. As a consequence, the set $\psi\left(B_{\sqrt{\rho}}\right)$ is bounded and, in view of $(a)$, non-convex.
Hence, each assumption of Theorem 2.7 is satisfied. Now, consider the operator $T:E\to E^*$ defined by
$$T(\alpha,\beta)(u,v)=\int_0^1\alpha(t)u(t)dt+\int_0^1\beta(t)v(t)dt$$
for all $(\alpha,\beta), (u,v)\in E$. Of course, $T$ is linear and the linear subspace $T(E)$ is total over $E$. Hence, $T(E)$ is
weakly-star dense in $E^*$.
 Moreover, notice that $T$ is continuous with respect to the weak-star topology of $E^*$. Indeed, let $\{(\alpha_n,\beta_n)\}$ be a sequence in $E$ converging to some $(\alpha,\beta)\in E$. Fix $(u,v)\in E$. We have to show that
$$\lim_{n\to \infty}T(\alpha_n,\beta_n)(u,v)=T(\alpha,\beta)(u,v)\ .\eqno{(14)}$$
Notice that
$$\lim_{n\to \infty}\left ( \int_0^1|\alpha_n(t)-\alpha(t)|dt + \int_0^1|\beta_n(t)-\beta(t)|dt\right ) = 0\ .\eqno{(15)}$$
On the other hand, we have
$$|T(\alpha_n,\beta_n)(u,v)-T(\alpha,\beta)(u,v)|=\left |\int_0^1(\alpha_n(t)-\alpha(t))u(t)dt+\int_0^1(\beta_n(t)-\beta(t))v(t)dt\right |$$
$$\leq \left ( \int_0^1|\alpha_n(t)-\alpha(t)|dt + \int_0^1|\beta_n(t)-\beta(t)|dt\right )\max\left\{\max_{[0,1]}|u|,\max_{[0,1]}|v|\right\}$$
and hence $(14)$ follows in view of $(15)$.
Finally, fix a convex set $S\subseteq C^0([0,1])\times C^0([0,1])$ dense in $C^0([0,1])\times C^0([0,1])$. Then, by the kind of continuity of
$T$ just now proved, the convex set
$T(-S)$ is weakly-star dense in $E^*$ and hence, thanks to Theorem 2.7, there exists $(\alpha_0,\beta_0)\in -S$ such that, if we put 
$$g(t,\xi)=\alpha_0(t)+\beta_0(t)f(\xi)\ ,$$
the functional $I-J_g$ has at least two critical points in $B_{\sqrt{\rho}}$ which are the claimed solutions of the problem in $(b)$, with 
$\alpha=-\alpha_0$ and $\beta=-\beta_0$.\par
Now, let us prove that $(b)\to (a)$. Assume that the restriction of $f$ to $\left [-{{\sqrt{\rho}}\over {2}},{{\sqrt{\rho}}\over {2}}\right ]$ is
constant. Let $c$ be such a value.  So, the classical solutions of the problem
$$\cases{-\omega\left(\int_0^1|u'(t)|^2dt\right)u''=c\beta(t)+\alpha(t) & in $[0,1]$\cr & \cr u(0)=u(1)=0\cr & \cr
\int_0^1|u'(t)|^2dt<\rho\cr}$$
are the critical points in $B_{\sqrt{\rho}}$ of the functional $u\to {{1}\over {2}}\tilde\omega\left(\int_0^1|u'(t)|^2dt\right)-
\int_0^1(c\alpha(t)+\beta(t))u(t)dt$. But, since $\omega$ is increasing and non-negative, this functional is strictly convex and so it possesses a unique critical point. The proof is complete.\hfill $\bigtriangleup$
\bigskip
\bigskip
{\bf Acknowledgements.} The author has been supported by the Gruppo Nazionale per l'Analisi Matematica, la Probabilit\`a e 
le loro Applicazioni (GNAMPA) of the Istituto Nazionale di Alta Matematica (INdAM) and by the Universit\`a degli Studi di Catania, PIACERI 2020-2022, Linea di intervento 2, Progetto ”MAFANE”. Thanks are also due to Prof. Lingju Kong whose pertinent questions lead to a clarification of the proof
of Theorem 2.10.
\vfill\eject
\centerline {\bf References}\par
\bigskip
\bigskip
\noindent
[1]\hskip 5pt A. R. ALIMOV and I. G. TSAR'KOV, {\it Connectedness and solarity in problems of best and near-best approximation}, Russian Math. Surveys, {\bf 71} (2016), 1-77.\par
\smallskip
\noindent
[2]\hskip 5pt V. S. BALAGANSK\u{I}I and L. P. VLASOV, {\it The problem of the convexity of Chebyshev sets}, Russian Math. Surveys, {\bf 51} (1996), 1127-1190.\par
\smallskip
\noindent
[3]\hskip 5pt N. V. EFIMOV and S. B. STE\u{C}KIN, {\it Approximative compactness and Chebyshev sets}, Dokl. Akad. Nauk SSSR, {\bf 140} (1961), 522-524.\par
\smallskip
\noindent
[4]\hskip 5pt F. FARACI and A. IANNIZZOTTO,  {\it An extension of a multiplicity theorem by Ricceri with an application to a class of quasilinear equations}, Studia Math., {\bf 172} (2006), 275-287.\par
\smallskip
\noindent
[5]\hskip 5pt F. FARACI and A. IANNIZZOTTO, {\it Well posed optimization problems and nonconvex Chebyshev sets in Hilbert spaces}, SIAM J. Optim., {\bf 19} (2008), 211-216.\par
\smallskip
\noindent
[6]\hskip 5pt  P. PUCCI and J. SERRIN, {\it A mountain pass theorem},
J. Differential Equations, {\bf 60} (1985), 142-149.\par
\smallskip
\noindent
\smallskip
\noindent
[7]\hskip 5pt B. RICCERI, {\it A general multiplicity theorem for certain nonlinear
equations in Hilbert spaces}, Proc. Amer. Math. Soc., {\bf 133} (2005),
3255-3261.\par
\smallskip
\noindent
[8]\hskip 5pt B. RICCERI, {\it A conjecture implying the existence of non-convex Chebyshev sets 
in infinite-dimensional Hilbert spaces},
Matematiche, {\bf 65} (2010), 193-199.\par
\smallskip
\noindent
[9]\hskip 5pt B. RICCERI, {\it On a minimax theorem: an improvement, a new proof and an overview of its applications},
Minimax Theory Appl., {\bf 2} (2017), 99-152.\par
\smallskip
\noindent
[10]\hskip 5pt I. G. TSAR'KOV, {\it Nonuniqueness of solutions of some differential equations and their connection with geometric approximation theory}, Math. Notes, {\bf 75} (2004), 259-271.\par
\smallskip
\noindent
[11]\hskip 5pt E. ZEIDLER, {\it Nonlinear functional analysis and its
applications}, vol. III, Springer-Verlag, 1985.\par
\bigskip
\bigskip
Department of Mathematics and Informatics\par
University of Catania\par
Viale A. Doria 6\par
95125 Catania, Italy\par
{\it e-mail address:} ricceri@dmi.unict.it

\bye